\documentclass{article}
\usepackage{amssymb,amsmath,
theorem,euscript}

\newcounter{sec}

\newcounter{punct}[sec]

\def\punct{\refstepcounter{punct}{\arabic{sec}.\arabic{punct}.  }}

\def\COUNTERS{\addtocounter{sec}{1}
              \setcounter{punct}{0}
          \setcounter{equation}{0}
          \setcounter{theorem}{0}
          \setcounter{problem}{0}
            \setcounter{Apunct}{0}
          }

\newtheorem{theorem}{Theorem}[sec]

\newtheorem{lemma}[theorem]{Lemma}

\newtheorem{observation}[theorem]{Observation}

\newtheorem{conjecture}{Conjecture}[sec]

\def\COUNTERS{\addtocounter{sec}{1}
              \setcounter{punct}{0}
          \setcounter{equation}{0}
          \setcounter{theorem}{0}
          }

\begin{document}

\def\SL{\mathrm {SL}}
\def\SU{\mathrm {SU}}
\def\GL{\mathrm  {GL}}
\def\U{\mathrm  U}
\def\OO{\mathrm  O}
\def\Sp{\mathrm  {Sp}}
\def\SO{\mathrm  {SO}}
\def\SOS{\mathrm {SO}^*}

\def\PGL{\mathrm  {PGL}}
\def\PU{\mathrm {PU}}

\def\Gr{\mathrm{Gr}}

\def\Fl{\mathrm{Fl}}

\def\OSp{\mathrm {OSp}}

\def\OU{\mathrm{OU}}

\def\OGL{\mathrm{OGL}}

\def\Mat{\mathrm{Mat}}

\def\Pfaff{\mathrm {Pfaff}}

\def\Ind{\mathrm{Ind}}

\def\B{\mathbf B}

\def\phi{\varphi}
\def\epsilon{\varepsilon}
\def\kappa{\varkappa}

\def\le{\leqslant}
\def\ge{\geqslant}

\renewcommand{\Re}{\mathop{\rm Re}\nolimits}

\renewcommand{\Im}{\mathop{\rm Im}\nolimits}

\def\pia{\pi_\downarrow}

\newcommand{\im}{\mathop{\rm im}\nolimits}
\newcommand{\indef}{\mathop{\rm indef}\nolimits}
\newcommand{\dom}{\mathop{\rm dom}\nolimits}
\newcommand{\codim}{\mathop{\rm codim}\nolimits}

\def\cA{\mathcal A}
\def\cB{\mathcal B}
\def\cC{\mathcal C}
\def\cD{\mathcal D}
\def\cE{\mathcal E}
\def\cF{\mathcal F}
\def\cG{\mathcal G}
\def\cH{\mathcal H}
\def\cJ{\mathcal J}
\def\cI{\mathcal I}
\def\cK{\mathcal K}
\def\cL{\mathcal L}
\def\cM{\mathcal M}
\def\cN{\mathcal N}
\def\cO{\mathcal O}
\def\cP{\mathcal P}
\def\cQ{\mathcal Q}
\def\cR{\mathcal R}
\def\cS{\mathcal S}
\def\cT{\mathcal T}
\def\cU{\mathcal U}
\def\cV{\mathcal V}
\def\cW{\mathcal W}
\def\cX{\mathcal X}
\def\cY{\mathcal Y}
\def\cZ{\mathcal Z}

\def\frA{\mathfrak A}
\def\frB{\mathfrak B}
\def\frC{\mathfrak C}
\def\frD{\mathfrak D}
\def\frE{\mathfrak E}
\def\frF{\mathfrak F}
\def\frG{\mathfrak G}
\def\frH{\mathfrak H}
\def\frJ{\mathfrak J}
\def\frK{\mathfrak K}
\def\frL{\mathfrak L}
\def\frM{\mathfrak M}
\def\frN{\mathfrak N}
\def\frO{\mathfrak O}
\def\frP{\mathfrak P}
\def\frQ{\mathfrak Q}
\def\frR{\mathfrak R}
\def\frS{\mathfrak S}
\def\frT{\mathfrak T}
\def\frU{\mathfrak U}
\def\frV{\mathfrak V}
\def\frW{\mathfrak W}
\def\frX{\mathfrak X}
\def\frY{\mathfrak Y}
\def\frZ{\mathfrak Z}

\def\fra{\mathfrak a}
\def\frb{\mathfrak b}
\def\frc{\mathfrak c}
\def\frd{\mathfrak d}
\def\fre{\mathfrak e}
\def\frf{\mathfrak f}
\def\frg{\mathfrak g}
\def\frh{\mathfrak h}
\def\fri{\mathfrak i}
\def\frj{\mathfrak j}
\def\frk{\mathfrak k}
\def\frl{\mathfrak l}
\def\frm{\mathfrak m}
\def\frn{\mathfrak n}
\def\fro{\mathfrak o}
\def\frp{\mathfrak p}
\def\frq{\mathfrak q}
\def\frr{\mathfrak r}
\def\frs{\mathfrak s}
\def\frt{\mathfrak t}
\def\fru{\mathfrak u}
\def\frv{\mathfrak v}
\def\frw{\mathfrak w}
\def\frx{\mathfrak x}
\def\fry{\mathfrak y}
\def\frz{\mathfrak z}

\def\fros{\mathfrak{s}}

\def\bfa{\mathbf a}
\def\bfb{\mathbf b}
\def\bfc{\mathbf c}
\def\bfd{\mathbf d}
\def\bfe{\mathbf e}
\def\bff{\mathbf f}
\def\bfg{\mathbf g}
\def\bfh{\mathbf h}
\def\bfi{\mathbf i}
\def\bfj{\mathbf j}
\def\bfk{\mathbf k}
\def\bfl{\mathbf l}
\def\bfm{\mathbf m}
\def\bfn{\mathbf n}
\def\bfo{\mathbf o}
\def\bfp{\mathbf q}
\def\bfr{\mathbf r}
\def\bfs{\mathbf s}
\def\bft{\mathbf t}
\def\bfu{\mathbf u}
\def\bfv{\mathbf v}
\def\bfw{\mathbf w}
\def\bfx{\mathbf x}
\def\bfy{\mathbf y}
\def\bfz{\mathbf z}

\def\bfA{\mathbf A}
\def\bfB{\mathbf B}
\def\bfC{\mathbf C}
\def\bfD{\mathbf D}
\def\bfE{\mathbf E}
\def\bfF{\mathbf F}
\def\bfG{\mathbf G}
\def\bfH{\mathbf H}
\def\bfI{\mathbf I}
\def\bfJ{\mathbf J}
\def\bfK{\mathbf K}
\def\bfL{\mathbf L}
\def\bfM{\mathbf M}
\def\bfN{\mathbf N}
\def\bfO{\mathbf O}
\def\bfP{\mathbf P}
\def\bfQ{\mathbf Q}
\def\bfR{\mathbf R}
\def\bfS{\mathbf S}
\def\bfT{\mathbf T}
\def\bfU{\mathbf U}
\def\bfV{\mathbf V}
\def\bfW{\mathbf W}
\def\bfX{\mathbf X}
\def\bfY{\mathbf Y}
\def\bfZ{\mathbf Z}

\def\R {{\mathbb R }}
 \def\C {{\mathbb C }}
  \def\Z{{\mathbb Z}}
  \def\H{{\mathbb H}}
\def\K{{\mathbb K}}
\def\N{{\mathbb N}}
\def\Q{{\mathbb Q}}
\def\A{{\mathbb A}}

\def\T{\mathbb T}

\def\bbA{\mathbb A}
\def\bbB{\mathbb B}
\def\bbD{\mathbb D}
\def\bbE{\mathbb E}
\def\bbF{\mathbb F}
\def\bbG{\mathbb G}
\def\bbI{\mathbb I}
\def\bbJ{\mathbb J}
\def\bbL{\mathbb L}
\def\bbM{\mathbb M}
\def\bbN{\mathbb N}
\def\bbO{\mathbb O}
\def\bbP{\mathbb P}
\def\bbQ{\mathbb Q}
\def\bbS{\mathbb S}
\def\bbT{\mathbb T}
\def\bbU{\mathbb U}
\def\bbV{\mathbb V}
\def\bbW{\mathbb W}
\def\bbX{\mathbb X}
\def\bbY{\mathbb Y}

 \def\ov{\overline}
\def\wt{\widetilde}
\def\wh{\widehat}

\def\P{\mathbb P}

\def\bO{\bf O}

\def\arr{\rightrightarrows}

\def\SS{\smallskip}

\def\ev{{\mathrm{even}}}
\def\od{{\mathrm{odd}}}

\def\q{\quad}

\def\F{\mathbf F}

\def\b{\mathbf b}

\def\RA{\Longrightarrow}

\begin{center}

\bf\Large Branching integrals and Casselman phenomenon

\sc\large

\bigskip

Yuri A. Neretin%
\footnote{Supported by the grant FWF, project P19064,
 Russian Federal Agency for Nuclear Energy,
Dutch  grant NWO.047.017.015, and grant JSPS-RFBR-07.01.91209 }

\end{center}

\bigskip

\begin{flushright}
To Mark Iosifovich Graev in his 85 birthday
\end{flushright}

{\small Let $G$ be a real semisimple Lie group,
$K$ its maximal complex subgroup,
and $G_\C$  its complexification.
  It is known that all the $K$-finite matrix
 elements   on $G$ admit holomorphic continuation
  to branching functions
  on $G_\C$ having singularities at the a prescribed
  divisor. We propose a geometric explanation
  of this phenomenon.}

\section{Introduction}

\COUNTERS

{\bf \punct Casselman theorem.} Let
$G$ be a real semisimple Lie group, let
$K$ be the maximal compact subgroup.
Let $G_\C$ be the complexification of $G$.

 Let $\rho$ be an infinite-dimensional
irreducible representation of $G$ in a
complete separable locally convex
  space $W$
\footnote{the case of unitary
representations in Hilbert spaces
is sufficiently non-trivial.}.
 Recall that a vector
$w\in W$ is $K$-finite if the orbit $\rho(G)v$
spans a finite
 dimensional subspace in
$W$.%
\footnote{Let us rephrase the definition.
We restrict $\rho$ to the subgroup $K$
and decompose the restriction into a direct
sum  $\sum V_i$  of finite-dimensional representations
of $K$. Finite sums of the form $\sum_{v_j\in V_j} v_j$
are precisely all the $K$-finite vectors.}

A $K$-finite matrix element is
 a function on $G$ of the form
$$
f(g)=\ell(\rho(g) v)
$$
where $v$ is a $K$-finite vector in
$W$ and
 $\ell$ is a $K$-finite linear functional,
 i.e., a $K$-finite element of the  dual representation.

\begin{theorem}%
\footnote{ Theorem 
 was obtained in famous preprints
 of W.Casselman on the Subrepresentation Theorem.
 Unfortunately, these works
  are unavailable for author;
 however they are included to the paper
 of W.Casselman and Dr.Milicic \cite{CM}.
There are (at least) two known proofs; the original
proof is based on properties
of system of partial differential
equations for matrix elements \cite{CM},
also the theorem can be reduced to
properties of Heckman--Opdam hypergeometric
functions \cite{HO} by a  simple trick \cite{Ner-harish}.}
\label{th:cass}
There is an (explicit) complex submanifold
$\Delta\subset G_\C$ of codimension 1 such that each
$K$-finite matrix element of $G$ admits a
continuation to an analytic multi-valued branching function
on $G_\C\setminus\Delta$.
\end{theorem}

{\sc Example.} Let $G=\SL(2,\R)$ be the group
of real matrices $\begin{pmatrix}a&b\\c&d\end{pmatrix}$,
whose determinant $=1$.
Then $K=\SO(2)$ consists of matrices
$
\begin{pmatrix}
\cos\phi&\sin\phi\\ -\sin\phi&\cos\phi
\end{pmatrix}
$, where $\phi\in\R$; the group $G_\C$
is the group of complex $2\times 2$ matrices
with determinant $=1$.
 The submanifold
$\Delta\subset\SL(2,\C)$ is a union of the
following four manifolds
\begin{equation}
a=0,\quad b=0,\quad c=0,\quad d=0
\label{eq:abcd}
\end{equation}
Indeed, in this case, there exists a canonical
$K$-eigenbasis.     All the matrix elements
in this basis
are Gauss hypergeometric
functions of the form
$$
{}_2F_1(\alpha,\beta;\gamma;\theta)
,\quad\text{where $\theta=\frac{ad}{bc}$}
$$
where the indices $\alpha$, $\beta$, $\gamma$
depend on parameters of a representation and
numbers of basis elements (see \cite{Vil}).

Branching points of ${}_2F_1$ are $\theta=0,1,\infty$.
Since $ad-bc=1$, only $\theta=0$ and $\theta=\infty$
are admissible; this implies (\ref{eq:abcd}).
\hfill$\square$

\medskip

Thus a representation $\rho$
of a real semisimple group admits
a continuation to an analytic matrix-valued function on
$G_\C$ having singularities at $\Delta$.
This fact seems to be strange if we look to
explicit constructions of representations.

Our purpose  is to clarify this phenomenon
and to find a direct geometric
construction of the analytic continuation.
We achieve this aim for a certain special case
(namely, for principal maximally degenerate
series of $\SL(n,\R)$, see Section 2)
and formulate a general conjecture
(Section 3). It seems that our explanation
(a reduction to the  'Thom isotopy Theorem'),
see \cite{Pha}, \cite{Vas})
is trivial. However, as far as I know it
is not known for experts in
the representation theory.

Section 4 is informal and contains a brief exposition
of various phenomena related to analytic continuations
of representations. Addendum contains a general discussion
of holomorphic continuations of
representations.


\SS

{\bf Acknowledgments.} The author is grateful
to Alexei Rosly and Victor Vasiliev for discussion
of this subject.


\section{Isotopy of cycles}

\COUNTERS

{\bf\punct  Principal degenerate series
for groups $\SL(n,\R)$.}
Let $G=\SL(n,\R)$ be the group of all
real  matrices with determinant $=1$.
 The maximal compact subgroup
  $K=\SO(n)$ is the group of all real orthogonal
 matrices.

Denote by $\R\P^{n-1}\subset\C\P^{n-1}$
the real and complex projective spaces;
recall that  the manifold $\R\P^{n-1}$
is orientable iff $n$ is even.

 Denote by $d\omega$ the $\SO(n)$-invariant Lebesgue measure
 on $\R\P^{n-1}$, let $d(\omega g)$ be its pushforward
 under the map $g$, denote by
 $$
J(g,x):=\frac{d\omega g}{d\omega}
 $$
the Jacobian of a transformation $g$ at a point
$x$.

Fix $\alpha\in\C$. Define a representation $T_\alpha(g)$
of the group $\SL(n,\R)$ in the space  $C^\infty(\R\P^{n-1})$
by the formula
$$
T_\alpha(g)f(x)=f(xg)J(g,x)^\alpha
$$
The representations $T_\alpha$ are called {\it representations
of principal degenerate series}.
If $\alpha\in \frac12+i\R$, then this representation
is unitary in $L^2(\R\P^{n-1})$.

\SS


{\bf\punct Discriminant submanifold $\Delta$.}
Denote by $g^t$ the transpose of a matrix $g$.
Denote by $\Delta$ the submanifold in $\SL(n,\C)$
consisting of matrices $g$ such that
the equation
$$
\det(gg^t-\lambda)=0
$$
has a multiple root.

\SS

{\it We wish to construct a continuation
of the function $g\mapsto T_\alpha(g)$
to a multi-valued function on $\SL(n,\C)\setminus\Delta$.}

\SS

For simplicity, {\it we assume $n$ is even.}%
\footnote{If $n$ is odd, then we must
replace the integrand in (\ref{eq:m-e}) by
a form on two sheet covering
of $\C\P^{n-1}\setminus\R\P^{n-1}$. Also we must replace
the cycle $\R\P^{n-1}$ by its two-sheet covering.}


\SS

{\bf\punct Invariant measure.}
Denote by $x_1:x_2:\dots:x_n$ the homogeneous
coordinates in a projective space.
The $\SO(n)$-invariant $(n-1)$-form
on $\R\P^{n-1}$ is given by
$$
d\omega(x)=(\sum_j x_j^2)^{-n/2}
\sum_j(-1)^{j} x_j\,dx_1\dots\wh{dx_j}\dots\, dx_n
$$
This expression can be regarded as a
meromorphic $(n-1)$-form on $\C\P^{n-1}$
having a pole on the quadric
$$
Q(x):=\sum x_j^2=0
$$

Now we can treat the Jacobian
$J(g,x)$ as a meromorphic function
on $\C\P^{n-1}$ having a zero at
the quadric $Q(x)=0$ and
a pole on the shifted quadric $Q(gx)=0$.

\SS


{\bf\punct $K$-finite functions.}
The following functions
span the space of $K$-finite functions on $\R\P^{n-1}$:
$$
f(x)=\frac{\prod x_j^{k_j}}
          {(\sum x_j^2)^{\sum k_j/2}},
          \qquad\text{where $\sum k_j$ is even}
$$
Evidently, they have singularities at
the quadric $Q(x)=0$ mentioned above.

\SS


{\bf\punct $K$-finite matrix elements.}
$K$-finite matrix elements are given by
formula
\begin{equation}
\{f_1,f_2\}=
\int_{\R\P^{n-1}}
f_1(x)f_2(xg) J(g,x)^\alpha d\omega(x)
\label{eq:m-e}
\end{equation}
The integrand is a holomorphic form on $\C\P^{n-1}$
 of the maximal
degree ramified over quadrics $Q(x)=0$,
$Q(xg)=0$. Denote by $\frU=\frU[g]$ the complement to these
quadrics.
Therefore locally in $\frU$
the integrand is a closed $(n-1)$-form.
Hence we can replace $\R\P^{n-1}$ by an arbitrary
isotopic cycle $C$ in $\frU$.

\SS


{\bf\punct Reduction to Pham Theorem.}
Now let $g(s)$ be a path in $\SL(n,\C)$ starting
in $\SL(n,\R)$. For each $s$ one has a pair $Q(x)=0$,
$Q(x\cdot g(s))=0$  of quadrics and the corresponding
complement $\frU(g(s))$.

\SS

 {\it Is it possible to construct
an isotopy $C(s)$ of the cycle $\R\P^{n-1}$
such that $C(s)\subset \frU(g(s))$ for all $s$?}

\SS

Now recall the following Pham theorem
(see F.Pham \cite{Pha}),  V.A.Vasiliev \cite{Vas}).

\begin{theorem}
Let $R_1(s)$, \dots, $R_l(s)$
be nonsingular complex hypersurfaces in $\C\P^k$
depending on a parameter. Assume that $R_j$ are transversal
(at all points for all values of the parameter $s$).
Then each cycle in the complement
to $\cup R_j(s)$ admits an isotopy according the parameter.
\end{theorem}


{\bf\punct Transversality of quadrics.}

\begin{lemma}
Let $A$, $B$ be non-degenerate symmetric matrices.
Assume that all the roots of the
characteristic  equation
$$
\det(A-\lambda B)=0
$$
are pairwise distinct.
Then
quadrics $\sum a_{ij}x_i x_j=0$ and $\sum b_{ij} x_i x_j=0$
are transversal.
\end{lemma}

By the Weierstrass theorem such pair of quadrics
can be reduced to
\begin{equation}
\sum \lambda_j x_j^2=0,\qquad \sum x_j^2=0
\label{eq:2}
\end{equation}
where $\lambda_j$ are the roots of the characteristic equation.
If they are not transversal at a point $x$, then
rank of the matrix
$$
\begin{pmatrix}
\lambda_1 x_1&\dots &\lambda_n x_n\\
x_1&\dots & x_n
\end{pmatrix}
$$
is 1. Therefore
\begin{equation}
(\lambda_i-\lambda_j)x_i x_j=0\quad\text{for all $i$, $j$}
\label{eq:1}
\end{equation}
The system (\ref{eq:1}), (\ref{eq:2}) is inconsistent.
\hfill$\square$


\SS

{\bf \punct Last step.} In our case, the matrices
of quadratic forms are $gg^t$ and $1$.
Therefore, by the virtue of the Pham Theorem a desired isotopy
of the cycle $\R\P^{n-1}$ exists.


\section{General case}

\COUNTERS


By the Subrepresentation Theorem, all the
irreducible representations of
a semisimple group $G$ are
subrepresentations of the  principal
(generally, non-unitary) series.
Therefore, it suffices to construct analytic
continuations for representations
 of the principal series.

For definiteness, we discuss the spherical
principal series
of
 the group $G=\SL(n,\R)$.

\SS


{\bf\punct Spherical principal series for $G=\SL(n,\R)$.}
Denote by $\Fl(\R^n)$ the space of all complete
flags of subspaces
$$
\cW:\,0\subset W_1\subset\dots\subset W_{n-1}\subset \R^n
$$
in $\R^n$; here $\dim W_k=k$.
By $\Gr_k(\R^n)$ we denote the Grassmannian
of all $k$-dimensional subspaces in $\R^n$.
By $\gamma_k$ we denote the natural
projection $\Fl(\R^n)\to\Gr_k(\R^n)$.

By $\omega_k$ we denote the
$\SO(n)$-invariant measure on $\Gr_k(\R^n)$.
For $g\in\GL(n,\R)$ we denote by
$J_k(g,V)$ the Jacobian of the transformation
$V\mapsto Vg$ of $\Gr_k(\R^n)$,
$$
J_k(g,V)=\frac{d\omega_k(Vg)}{d\omega_k(V)}
$$

Fix $\alpha_1$, \dots, $\alpha_{n-1}\in\C$.
The representation $T_\alpha$
of the spherical principal series
of the group $\SL(n,\R)$ acts in the space
$C^\infty(\Fl(\R^n))$ by the formula
$$
T_\alpha(g)f(\cW)=f(\cW\cdot g)
\prod_{k=1}^{n-1} J_k(g,\gamma_k(\cW))^{\alpha_k}
$$


\SS

{\bf\punct Singularities.}
Consider the symmetric bilinear form in
$\C^n$ given by
$$
B(x,y)=\sum x_j y_j
$$
By $L_k\subset \Gr_k(\C^n)$ we denote the set of
all the $k$-dimensional subspaces, where the form
$B$ is degenerate%
\footnote{Equivalently, we can consider all the
$(k-1)$-dimensional subspaces in $\C\P^{n-1}$
tangent to the quadric $\sum x_j^2=0$.}.
 By $\cL\subset \Fl(\C^n)$
we denote the set of all the flags
$W_1\subset\dots\subset W_{n-1}$, where
$W_k\in L_k$ for some $k$.

\SS

In fact,  all the $K$-finite
functions on $\Fl(\R^n)$ admit  analytic
 continuations to $\Fl(\C^n)\setminus \cL$
(a singularity on $\cL$ is a pole or two-sheet branching).

\SS


{\bf\punct A conjecture.}

\begin{conjecture}
Let $\gamma(t)$ be
 a path on $\GL(n,\C)$
avoiding the discriminant submanifold $\Delta$,
let $\gamma(0)\in\SL(n,\R)$.
Then there is an isotopy $C(t)$ of the cycle
$\Fl(\R^n)$ in the space
$\Fl(\C^n)$ avoiding the submanifolds
$\cL$ and $\cL\cdot g(s)$
\end{conjecture}

Such isotopy  produces an analytic continuation
of representations of principal series of $\SL(n,\R)$.


\vspace{30pt}

{\bf \Large Addendum. Survey of holomorphic continuations
of representations}

\bigskip

{\it Let $G$ be a connected linear Lie group}, i.e., a Lie group
admitting an embedding to a matrix group.
Denote by $G_\C$ its complexification.
Let $\rho$ be an irreducible representation
of $G$ (in a Frechet space). We are interested in the following problems:

\SS

--- Is it possible to extend $\rho$ holomorphically to $G_\C$?

\SS

--- Is it possible to extend $\rho$ holomorphically to an
open domain $U\subset G_\C$.

\SS

See, also, \cite{Ner-book}, Section 1.5.


\SS

{\bf A.1. Weyl trick.} Let $\rho$ be a finite dimensional
representation of a semisimple Lie group $G$.
Then $\rho$ admits the holomorphic continuation to
the group $G_\C$.

\SS


{\bf A.2. Why the Weyl trick does not survive for infinite-dimensional unitary
representations?}  Let $G$ be a noncompact Lie group, let $\rho$ be its 
 irreducible faithful unitary
representation. Let $X$ be a noncentral element of the Lie algebra $\frg$.
It is more-or-less obvious that
 the operator $\rho(X)$ is  unbounded%
\footnote{I propose a collection of arguments that can be combined to a proof
in various way. First, this is valid for the two-dimensional non-Abelian algebra
(and therefore this is valid for elements of $\frg$ that can be included to such
a subalgebra).
Second, this is valid for nilpotent algebras, we apply  Kirillov's monomial induction theorem, see \cite{Kir}, then the Lie algebra acts by first 
order differential operators. which are obviously unbounded. Third,
for we know that for semisimple Lie groups $K$-spectra are always unbounded.
 Fourth, elements $X$ with bounded $\rho(X)$ form an ideal. 
In fifth, we can consider its normalizer}.

Then, for $t$, $s\in\R$,
$$
\rho(\exp(t+is)X)=\exp(isX) \exp(tX)
$$
Since $iX$ is self-adjoint, then
$\exp(tX)$ is unitary; on the other hand $\exp(isX)$
have to be unbounded for all positive $s$ or for all
negative $s$(and usually it is unbounded for all $s$..

However, this argument does not removes completely
an idea of holomorphic continuation,  since it remain two following
logical possibilities

\SS

--- a holomorphic extension exists in spite of the unboundedness of operators.

\SS

--- If a specter of $X$ is contained on a positive half-line,
then $\exp(tX)$ is defined for negative $t$. We can hope to construct
something from elements of this kind.

\SS

The second variant is realized for Olshanski semigroups,
see below,
the first variant is general, this follows from the Nelson Theorem.


\SS

{\bf A.3. Nelson's paper.} In 1959 E.Nelson \cite{Nel}
 proved  that each
unitary irreducible representation $\rho$
of a real Lie group $G$ has a dense set of analytic vectors.
This implies that $\rho$
can be extended analytically
to a sufficiently small neighborhood of $G$ in $G_\C$.

However, usually this continuation can be done in a constructive
way as it is explained below (see also \cite{Ner-book}, Section 1.5.).

\SS


{\bf A.4. Induced representations.} 
First, we recall a definition of induced representations.

Consider a  Lie group $G$ and its closed
connected subgroup $H$. Let $\rho$ be a
representation of $H$ in a {\it finite-dimensional}
complex space $V$.

These data allow to construct canonically a
 vector bundle ({\it skew product}) over $G/H$ with a fiber $V$.
Recall a construction (see, for instance, \cite{Kir}).
 Consider the direct product
$G\times V$ and the equivalence relation
$$(g,v)\sim (gh^{-1},\rho(h)v),
\qquad\text{where $g\in G$, $v\in V$, $h\in H$}
$$
Denote by $R=G\times_V H$  the quotient space. The standard map
$G\to G/H$ determines a map $R\to G/H$
(we simply forget $v$). A fiber can be (noncanonically) identified
with $V$.

Next, the group $G$ acts on $G\times V$ by transformations
$$(g,v)\mapsto (rg,v),\qquad
\text{where $g\in G$, $v\in V$, $r\in G$}
$$
This action induces the action of $G$ on $G\times_H V$,
the last action commutes with the projection
$G\times_H V\to G/H$.

Therefore $G$ acts in the space of sections of the bundle
$G\times_H V\to G/H$ (because graph of a section is a
subset in the total space; the group $G$ simply move subsets).
{\it Induced representation $\pi=\Ind_H^G(\rho)$
 is the representation
of $G$ in a space of sections of the bundle
 $G\times_H V\to G/H$.}

 The most important example are principal series,
 which were partially discussed above.

Our definition is not  satisfactory since rather often
 it is necessary to specify the space of sections
(for instance, smooth functions,
$L^2$-functions, distributions, etc.). This discussion
is far beyond our purposes, for the moment
let us consider the space $C^\infty[G/H;\rho]$ of smooth sections.

\SS


{\bf A.5. Analytic continuation of induced representations.}
Here we discuss some heuristic
arguments (see \cite{Ner-book}, Section 1.5).
 Their actual usage
depends on explicit situation.

Denote by $H_{[\C]}$ the complexification
of the group $H$ inside $G_\C$.%
\footnote{For elements $X$ of the Lie algebra of $H$ we consider
the subgroup in $G_\C$ spanned by all $\exp((t+is)X)$.}

The (finite-dimensional) representation $\rho$
admits a holomorphic extension to a representation
of the universal covering group
of $\wt H_{[\C]}$ of  $H_{[\C]}$ in the  space $V$.
For a moment, let us require two assumptions%
\footnote{The second assumption is very restrictive. It is not
hold for the parabolic induction.}

\SS

--- {\it  $H_{[\C]}$ is closed in $G_\C$.} 

\SS

--- {$\rho$ is a linear representation
of $ H_{[\C]}$.}

\SS

Under these assumptions, the same 
construction of skew-product produces the
 bundle $(G_\C)\times_{H_\C} V\to G_\C/ H_\C$. Moreover,
$$
G\times_H V\,\subset\,(G_\C)\times_{H_\C} V
$$

Now let us agree with the next assumption%
\footnote{If it is not hold, then we go to  subsection A.7,
where all the assumptions are omitted.}.
{\it Let the space of holomorphic sections of
$(G_\C)\times_{H_\C} V$ be dense in $C^\infty$ on $G/H$.}
 Then
we get a  holomorphic continuation of the induced
representation $\pi$ to the whole complex group
 $G_\C$. More precisely, we slightly reduce
the space of representation, but 'formulae' determining
a representation are  the same.

\SS

This variant is realized for all the nilpotent Lie groups.

\SS


{\bf A.6. Nilpotent Lie groups.}

\SS

{\sc Example.} Let $a$, $b$, $c$ range in $\R$, 
Consider operators 
$T(a,b,c)$ 
in $L^2(\R)$ given by
$$
T(a,b,c)f(x)=f(x+a)e^{ibx+c}
$$
They form a 3-dimensional group, namely the Heisenberg group.
Now let $a$, $b$, $c$ range in $\C$ and $f$ ranges in the space
$Hol(\C)$ of
 entire functions. Then the
same formula determines a representation of the complex Heisenberg group
in $Hol(\C)$. After this operation, the space of representation was completely changed.
However it is easy to invent a dense subspace in $L^2(\R)$ consisting
of holomorphic functions and invariant with respect to all the operators
$T(a,b,c)$.   \hfill $\square$

\SS

Now, let $G$ be a nilpotent Lie group.
By the Kirillov Theorem \cite{Kir-orbit},
each unitary representation of a nilpotent Lie group
$G$ is induced from an one-dimensional representation of
a subgroup $H$; the  manifolds $G/H$ are equivalent
to standard spaces $\R^m$. Therefore, $G_\C/H_\C$
is the standard complex space $\C^m$, and we get a representation
of $G_\C$ in the space of entire functions.

However, the following R.Goodman--G.L.Litvinov Theorem
(R.Goodman \cite{Goo}, G.L.Litvinov
\cite{Lit1}, \cite{Lit2})
is more delicate.

\SS

{\bf Theorem.} {\it
Let $\rho$ be an irreducible unitary representation 
of a nilpotent group $G$ in a Hilbert space $W$.
There exists a (noncanonical) dense subspace $Y$ with its own 
Frechet topology
and holomorphic representation
$\wt \rho$  of $G$ in the space $W$ 
coinciding with $\rho$ on $G$.}

\SS

Let us explain how to produce a subspace $Y$.
Let $\rho$ be a unitary representation of a nilpotent
group $G$ in a space $W$.
Let $f$ be an entire function on $G_\C$ (it is specified below).
 Consider the operator
$$
\rho(f)=\int_G f(g) \rho(g) \,dg
$$
Let $r\in G_\C$. We write formally
$$
\rho(r)\rho(f)=\int_G f(g)\rho(r g)\,dg
=\int_G f(gr^{-1})\rho(g)
$$
{\it Assume that for each $r\in G_\C$ the function
$\gamma_r(g):=f(gr^{-1})$ is integrable on $G$.}
Under this condition we can define
operators
$$
\rho(r):\,\Bigl\{\text{Image of $\rho(f)$}\Bigr\}
\,\to\, W
$$
as just now.

In fact, we need  a subspace $Z$ in $L^1(G)$ 
consisting of entire functions and invariant with respect to
complex shifts. To be sure that 
$$
Y:=\cup_{f\in Z}
\mathrm{Im}(\rho(f))\subset W
$$
is dense, we need in a sequence of positive $f_j\in Z $ converging
to $\delta$-functions; then $\rho(f_j)v$ converges to $v$ for all $v\in W$.
In what follows we  describe a simple trick that allows to construct many functions
$f$ and a subspace $Z$.

First, let $G=T_n$ be the unipotent upper triangular subgroup
of order $n$. Let $t_{ij}$, $i<j$,
be the natural coordinates on $T_n$. Write them in the order
$$
t_{12},t_{23}, t_{34},\dots, t_{(n-1)n},
t_{13},t_{24},\dots, t_{(n-2)n},t_{14},\dots
$$
and re-denote these coordinates by $x_1$, $x_2$, $x_3$, \dots.
In this coordinates, the right shift $g\mapsto gr^{-1}$
is an affine transformation of the form
$$
(x_1,x_2,x_3,\dots)\mapsto (x_1+a_1, x_2+a_2+b_{21}x_1,
x_3+a_3+b_{31}x_1+b_{32}x_2,\dots)
$$
Now we can choose a desired function $f$ in the form
$$
f(x)=\exp\Bigl\{-\sum p_j x_j^{ 2 k_j}\Bigr\},
\quad \text{where $p_j>0$ and $k_1>k_2>\dots>0$ are integers}
$$

For an arbitrary nilpotent $G$ we apply the Ado theorem
(in fact, the standard proof, see \cite{Ner-ado})
 produces a polynomial embedding
of $G$ to  some $T_n$ with very large $n$).

\SS


{\bf A.7. Local holomorphic continuations.} Now let $G\supset H$
be arbitrary.

The construction of a skew product
survives locally.  It determines a holomorphic bundle
on a (noncanonical)
 neighborhood $U\subset G_\C/H_\C$ of $G/H$.
 Denote by $\cA(U)$ the space of holomorphic sections
 of this bundle. Let $r\in G_\C$ satisfies
 $r\cdot G/H\subset U$. The $r$ induces an operator
 $\pi(r):\cA(U)\to C^\infty$ and we obtain an
 approach to local analytic continuation
 induced representation.

 \SS


 {\bf A.8. Local analytic continuations
 for semisimple groups.} For definiteness, consider
  $G=\SL(n,\R)$. Denote by $P$
 be the minimal parabolic (i.e., $P$ is the group
 of upper triangular matrices), Then $G/P$ is
the flag space mentioned above. Next, $G_\C=\SL(n,\C)$,
and $P_{[\C]}$ is the group of complex upper-triangular matrices
of the form
$$
B=
\begin{pmatrix} b_{11}&b_{12}&\dots\\
                 0 & b_{22}& \dots\\
                 \vdots&\vdots&\ddots
                 \end{pmatrix}
\text{$a_{ii}\ne 0$, other $a_{ij}$ are arbitrary.}
$$
Evidently, the group $P_\C$ is not simply connected.

Fix $s_j\in \C$ and consider the one-dimensional character
$$
\chi_s(B):=\prod_{j=1}^n b_{jj}^{s_j}
$$
of $P$. The function $\chi_s$
is defined on $P_\C$ only locally,
however this is sufficient
for the  arguments of the previous subsection.

This kind of arguments can be easily applied
to an arbitrary representation of a nondegenerate
principal series. Keeping in the mind Subquotient
Theorem, we easily get the following statement.

\begin{observation}
Let $G$ be a semisimple Lie group. Then
there are (noncanonical) open sets
$U_1\subset U_2\subset G_\C$
containing $G$ such that $U_2\supset U_1\cdot U_1$
and the following property holds.
Let $\rho$ be an irreducible representation of
$G$ in a Frechet space $W$. Then there is a (noncanonical)
dense subspace $Y\subset W$ (equipped with
its
own Frechet topology) and an operator-valued holomorphic
function $\wt \rho:U_2\to \mathrm{Hom}(Y,W)$ such that
$\wt\rho=\rho$ on $G$ and
$$
\wt\rho(g_1)\wt\rho(g_2)y=
\wt\rho(g_1g_2)y
\qquad\text{for $g_1$, $g_2\in U_1$, $y\in Y$}
$$
\end{observation}

Certainly, the operators $\wt\rho(g)$ are unbounded
in the topology of $\mathrm{Hom}(W,W)$.


\SS


{\bf A.9. Crown.}
 D.N.Akhiezer and S.G.Gindikin
(see \cite{AG}) constructed a certain explicit
domain $\cA\subset G_\C$ ('crown')
 to which all the spherical functions of a real semisimple
 group $G$ can be extended.
  Also the crown is a domain of holomorphy of all
 irreducible representations of  $G$, see
  B.Krotz, R.Stanton,
\cite{KS1}, \cite{KS2}.

Relation of their constructions with
our previous considerations are not completely clear.

\SS



{\bf A.10. Olshanski semigroup.} In all the previous
examples, operators of holomorphic continuation
are unbounded in the initial topology. There is an important
exception.

 Unitary highest weight representations of
a semisimple Lie group $G$ admit  holomorphic
continuations to a certain subsemigroup $\Gamma\subset G_\C$
(M.I.Graev \cite{Gra}, G.I.Olshanski \cite{Ols}).
Since this situation is well understood, we omit further
discussion, see also \cite{Ner-gauss}.

\SS


{\bf A.11. Infinite-dimensional groups.}
Induction (in different variants) is the main tool
of construction of representations of Lie groups.

For infinite-dimensional groups the induction exists%
\footnote{If a group acts on the space with measure, then
it acts in $L^2$.}
but it is a secondary
tool (however, the algebraic variant of
induction is important for infinite dimensional Lie algebras).
 A more effective instrument are symplectic and orthogonal
spinors.

Let us realize the standard real orthogonal group $\OO(2n)$
as a group 
of $(n+n)\times(n+n)$ complex matrices $g$ having the structure
$$
g=
\begin{pmatrix}
\Phi&\Psi\\
\ov\Psi&\ov\Phi
\end{pmatrix}
$$
that are orthogonal in the following sense
$$
g\begin{pmatrix}
0&1\\1&0 
\end{pmatrix}g^t
=
\begin{pmatrix}
0&1\\1&0 
\end{pmatrix}
$$
Actually we do the following. Consider the space $\R^{2n}$ equipped with a standard
basis $e_1$, \dots,  $e_n$, $e_{n+1}$,\dots, $e_{2n}$.
Then we pass to the space $\C^{2n}$ and write matrices
of real orthogonal operators in the basis 
$$
e_1+ie_{n+1}, \, e_2+ ie_{n+2},\,\dots,\, e_{n}+i e_{2n},
e_1-ie_{n+1}, \, e_2-ie_{n+2},\,\dots,\, e_{n}-i e_{2n}
$$

Next, we put $n=\infty$. Denote by $\OU(2\infty)$
the group of all the bounded matrices of the same structure
satisfying an additional condition:
$\Psi$ is a Hilbert--Schmidt matrix (i.e., the sum  $\sum |\psi_{kl}|^2$
of squares of matrix elements is finite).

By the well-known theorem of F.A.Berezin \cite{Ber}, \cite{Ber-book}
and D.Shale-W.Stinespring \cite{SS}, the spinor representation
is well-defined on the group $\OU(2\infty)$. 

Numerous infinite dimensional groups $G$
can be embedded in a natural way to $\OU(2\infty)$,
after this we
can restrict the spinor representation to $G$.
For instance, for the loop group $C^\infty(S^1,\SO(2n)]$
we consider the natural action in the space
$L^2(S^1,\R^{2n})$ and define the operator
of the complex structure in this space via Hilbert transform
(see, for instance, \cite{Ner-book}).
Applying the spinor representation, we get the so-called
basic representation of the loop group. For production
of other highest weight representations
 we apply restrictions and tensoring.

However, the group $\OU(2\infty)$ admits 
a complexification $\OGL(2\infty,\C)$ consisting
of complex matrices $g=\begin{pmatrix} A&B\\C&D\end{pmatrix}$
that are orthogonal in the same sense with Hilbert--Schmidt blocks
$B$ and $C$. The spinor representation of
$\OU(2\infty)$ admits a holomorphic continuation to the complex
group $\OGL(2\infty)$, see \cite{Ner-spinor}, \cite{Ner-book},
 the operators
of continuation are unbounded in the initial topology, but are bounded
on a certain dense Frechet subspace equipped with its own topology.

This produces highest weight representation
of complex loop groups as free byproducts
(see another approach in R.Goodman, N.Wallach \cite{GW}).

More interesting phenomenon arises for
the group $\mathrm{Diff}$ of diffeomorphisms of circle,
in this case the analytic continuation exists in spite
of nonexistence $\mathrm{Diff}_\C$, see \cite{Nerb}.

{\tt Math.Dept., University of Vienna,

 Nordbergstrasse, 15,
Vienna, Austria

\&

Institute for Theoretical and Experimental Physics,

Bolshaya Cheremushkinskaya, 25, Moscow 117259,
Russia

e-mail: neretin(at) mccme.ru

URL:www.mat.univie.ac.at/$\sim$neretin
}

\end{document}